\documentstyle[11pt]{article}
\newcommand{\dubbelR}{{\sf I}\kern-.12em{\sf R}}
\newcommand{\dubbelN}{{\sf I}\kern-.12em{\sf N}}
\newcommand{\dubbelZ}{\mbox{\sf Z\hspace*{-5.4pt}Z}}
\newfont{\Bbb}{msbm10 scaled\magstep1}
\newcommand{\CC}{\hbox{\Bbb C}}

\def\half{{\frac{1}{2}}}
\newcommand{\ij}{i\kern -0.08em j}

\newcommand{\be}{\begin{equation}}
\newcommand{\ee}{\end{equation}}
\newcommand{\bea}{\begin{eqnarray}}
\newcommand{\eea}{\end{eqnarray}} 
\newcommand{\nn}{\nonumber}
\newcommand{\dubbelint}{\mathop{\int\int}\limits_{\kern-5.5pt {V-V_1}}}
\newcommand{\al}{\alpha}

\newcommand{\dl}{\delta}
\newcommand{\Dl}{\Delta}
\newcommand{\ep}{\epsilon}

\newcommand{\lm}{\lambda}

\newcommand{\xinv}{x^{-1}}
\newcommand{\yinv}{y^{-1}}
\newcommand{\zinv}{z^{-1}}

\newcommand{\ninv}{n^{-1}}
\newcommand{\fox}{f \otimes x}
\newcommand{\age}{g_A}
\newcommand{\DA}{{\cal D}({\cal A})}
\newcommand{\DG}{{\cal D}(G)}
\newcommand{\HAal}{{\cal H}^A_{\al}}
\newcommand{\VAal}{V^A_{\al}}
\newcommand{\PAal}{\Pi^A_{\al}}

\newcommand{\oba}{\otimes_{{\cal B}_A}}
\newtheorem{thm}{Theorem}[section]
\newtheorem{lem}[thm]{Lemma}
\newtheorem{prop}[thm]{Proposition}
\newtheorem{cor}[thm]{Corollary}
\newtheorem{defn}[thm]{Definition}
\newtheorem{rem}[thm]{Remark}
\begin{document}
\title{
\hfill{\normalsize UvA-WINS-Wisk.\ 96-08}\\[-0.3 cm]
\hfill{\normalsize UvA-WINS-ITFA 96-19}\\[-0.3 cm]
\hfill{\normalsize q-alg/9605044}\\[-0.3 cm]
\hfill{\normalsize to appear in Journal of Lie Theory}\\[1.5 cm]
       The Quantum Double of a (locally) Compact Group}
\author{T.H. Koornwinder \thanks{email: thk@wins.uva.nl}
      \\Department of Mathematics, University of Amsterdam, \\
        Plantage Muidergracht 24,
        1018 TV Amsterdam, The Netherlands.
 \and N.M. Muller\thanks{email: nmuller@phys.uva.nl}
      \\Institute for Theoretical Physicis, University of Amsterdam, \\
        Valckenierstraat 65,
        1018 XE Amsterdam, The Netherlands.}
\date{9 September, 1996.}
\maketitle
\begin{abstract}
We generalise the quantum double construction of Drinfel'd  to the 
case of the (Hopf) algebra of suitable functions on a compact or
locally compact group. We will concentrate on the $*$-algebra structure
of the quantum double. If the conjugacy classes in the group are 
countably separated, then we classify the irreducible
$*$-representations by using the connection with so--called transformation
group algebras. For finite groups, we will compare our description to the
result of Dijkgraaf, Pasquier and Roche. Finally we will work 
out the explicit examples of $SU(2)$ and $SL(2,\dubbelR)$.  
\end{abstract}
\section*{Introduction}
The quantum double of a Hopf algebra (or, double quantum group) was
introduced by Drinfel'd in \cite{Drin}. Quantum doubles are important 
examples of quasitriangular
(quasi) Hopf algebras, and in that sense they are well-studied, see for 
instance  \cite{CharPres}, \cite{ShniSte}, \cite{Majbook}. 

The existing theory of quantum doubles has beautiful applications in
physics: in \cite{DPR} Dijkgraaf, Pasquier and Roche show that the 
representation theory covers the main interesting data  of particular
orbifolds of Rational Conformal Field Theories. Tightly connected are the
topological interactions in spontaneously broken gauge theories. In
\cite{BDW1}, \cite{BDW2}  Bais, Van Driel and De Wild Propitius show that
the non-trivial fusion and braiding properties of the excited states in
broken gauge theories can be fully described by the representation theory
of a quantum double which is constructed from a finite group $G$
via a finite dimensional Hopf algebra ${\cal A}=\CC[G]$ (the corresponding 
group ring). For a detailed treatment, see \cite{Mark}. 

{}From a physical as well as a mathematical point of view it is natural to 
ask whether the quantum double construction from a finite group can be
generalized to the case of a compact, or even locally compact group. 
In this report we will explicitly
construct the quantum double $\DG$ corresponding to a (locally) compact 
group $G$. It has a natural $*$-structure, and we will find a class of 
$*$-representations (unitary representations), and prove rigorously that 
they form a complete set of irreducible $*$-representations. The construction 
uses the representation theory of {\sl transformation group algebras}, which we
will discuss in detail, and its connections with the theory of induced
representations of locally compact groups via the imprimitivity theorem
of Mackey \cite{Mack49}. For an overview of the theory of (induced)
representations of locally compact groups, see for instance Chapters 
9--10--11 in \cite{Mee79}. In fact, the same construction works more generally
for classifying the irreducible $*$-representations of transformation group
algebras, see Glimm \cite{gli62}. 

To be more precise, the construction is done in the following way: For the
generalization of the quantum double, we choose the algebra $C_c(G\times G)$
of continuous functions on $G\times G$ with compact support. This allows us
to use the representation theory of transformation group algebras 
$C_c(X\times G)$, where the locally compact group $G$ acts continuously on a 
locally compact space $X$. Under the technical (but crucial) assumption
that the conjugate action of $G$ on $G$ is {\sl countably separated}, 
classification
of the irreducible $*$-representations of these algebras turns out to be 
equivalent to classifying the irreducible representations $(\tau,P)$ of the 
pair $(G,{\cal O})$, with ${\cal O}$ an orbit of $G$ on $X$. Writing 
$G/H\simeq {\cal O}$, with
$H$ a closed subgroup of $G$, then $(G, G/H, P)$ is a {\sl system of 
imprimitivity} for the unitary representation $\tau$ of $G$ on a Hilbert space
$V$, and $P$ is a projection valued measure on $G/H$ acting on $V$. From
Mackey's imprimitivity theorem it follows that such representations $\tau$
are precisely the representations of $G$ which are
induced from  unitary irreducible representations
$\al$ of $H$. The classification of irreducible $*$-representations of the
quantum double is a direct consequence.

We will show that the case of finite $G$ 
is covered by our description, and it leads to representations which are 
isomorphic to the ones derived in \cite{DPR}. 

Finally we will work out the explicit examples of $G=SU(2)$ (compact) and
$G=SL(2,\dubbelR)$ (non--compact). Their interesting applications in 
physics will be discussed in a forthcoming paper by the second author, 
where the connection with a quantization of a Chern--Simons theory in 
(2+1) dimensions will be described.

We are also studying in detail the coalgebra structure of the 
continuous quantum double, and in a follow--up of this paper we 
will give the tensor product decomposition (`fusion rules') for the
quantum double representations, the universal $R$-matrix, and its
action on a tensor product state. 

In fact, questions about fusion rules were our original motivation for
this work. However,
it soon became apparent that even the definition of the quantum double of
$C(G)$, and the classification of the irreducible representations of this
quantum double had to be clarified. This led us to a thorough study of
quite some older literature on transformation group algebras, which
turned out to be well applicable for our case of the quantum double. To our
knowledge these references have not been put together in this combination 
before.

We do not claim originality in the contents of our main results.
Notably, the main Theorem 3.9 occurs in Glimm [16], in a somewhat hidden way.
Our reformulation of the theorem may be more suitable for applications,
for instance in physics, and makes it possible to treat concrete examples.
For expository reasons we have added our own version of the proof, with
emphasis on the link with the imprimitivity theorem. It gives insight
into the way we have derived the characterisation and classification
of the irreducible unitary representations of the quantum double. The
latter is
necessary for the computation of the fusion rules and for the braiding
properties ($R$-matrix) of the model.
\section{Construction}
Drinfel'd \cite{Drin} gives the following definition of the 
{\sl quantum double}
$\DA$ of a Hopf algebra ${\cal A}$. (This definition is only 
mathematically precise if ${\cal A}$ is a finite dimensional Hopf algebra.
However, there is a way out by working with a dual pair of bialgebras,
see Majid \cite{Majbook}, p.296.)
\begin{defn}
Let ${\cal A}$ be a Hopf algebra and ${\cal A}^0$ 
the dual Hopf algebra to ${\cal A}$ with the opposite comultiplication. 
Then $\DA$ is the unique quasi-triangular Hopf algebra with universal
R-matrix $R\in \DA\otimes \DA$ such that
\begin{itemize}
\item[i.]  As a vector space, $\DA = {\cal A}\otimes{\cal A}^0$.
\item[ii.] ${\cal A} = {\cal A}\otimes 1$ and ${\cal A}^0 = 1 \otimes
{\cal A}^0$ are Hopf subalgebras of $\DA$.
\item[iii.] The mapping $x \otimes \xi \mapsto x\xi : {\cal A}\otimes
{\cal A}^0 \to \DA$ is an isomorphism of vector spaces. Here $x\xi$ is
short notation for the product $(x\otimes 1)(1\otimes\xi)$. 
\item[iv.] Let $(e_i)_{i\in I}$ be a basis of ${\cal A}$ and 
$(e^i)_{i\in I}$ the dual basis of ${\cal A}^0$. Then 
\be
R = \sum_{i\in I} (e_i \otimes 1)\otimes
(1\otimes e^i),
\ee
independent of the choice of the basis. 
\end{itemize}
\end{defn}
Tensor products are taken over the field $\CC$. Now write
\bea
\Dl (x) = \sum_{(x)} x_{(1)} \otimes x_{(2)} \nn \\
(\Dl \otimes id) \circ \Dl (x) = \sum_{(x)} x_{(1)} \otimes x_{(2)} \otimes
x_{(3)}
\eea
for (iterated) comultiplication on ${\cal A}$, and similarly for 
comultiplication on ${\cal A}^0$. Then the comultiplication on $\DA$ is 
given by
\be
\Dl (x\xi) = \sum_{(x),(\xi)} x_{(1)}\xi_{(1)} \otimes x_{(2)}\xi_{(2)}, \qquad
(x\in {\cal A}, \xi \in {\cal A}^0). 
\label{eq:comultda}
\ee
and the multiplication by
\be
\xi x =\sum_{(x),(\xi)} \xi_{(1)}(Sx_{(1)})\, \xi_{(3)}(x_{(3)})\,
x_{(2)}\xi_{(2)}, \qquad (x\in {\cal A}, \xi \in {\cal A}^0) 
\label{eq:multda}
\ee
where $S$ denotes the antipode. 

We consider the case where $G$ is a finite group and ${\cal A} := C(G)$, 
the space of all complex
valued functions on $G$, which becomes a Hopf algebra under pointwise 
multiplication, with comultiplication
\be
(\Dl f)(x,y) := f(xy), \qquad (x,y \in G) 
\ee
and with antipode
\be
(Sf)(x) := f(\xinv), \qquad (x \in G).
\ee
The dual of $C(G)$ is the group algebra $\CC[G]$ with comultiplication
\be
\Dl (x) = x \otimes x, \qquad (x\in G).
\ee
The pairing is given by
\be
\langle f,x \rangle = f(x), \qquad (f\in C(G), x \in G).
\ee
\vspace*{.2cm}

For the quantum double of $C(G)$ we now write
\be
{\cal D}(G) := {\cal D}(C(G)) = C(G) \otimes \CC[G] \simeq C(G,\CC[G]).
\ee
Thus $f\otimes x \in C(G) \otimes \CC[G]$ can be considered as the mapping 
\be
z \mapsto f(z)x : G \to \CC[G].
\ee
Also ${\cal D}(G) \otimes {\cal D}(G) \simeq C(G\times G, \CC[G]\otimes 
\CC[G])$. Now formulas (\ref{eq:comultda}) and (\ref{eq:multda}) can be
written as:
\bea
(\Dl (f\otimes x))(y,z) = f(yz) x\otimes x, \qquad (x,y,z \in G, f\in C(G))
\label{eq:comdub} \\
(1 \otimes x)(f \otimes e) = f(\xinv .\,x)\otimes x, \qquad (e\,\,
\mbox{unit in}\, G). 
\eea
Hence
\be
(f_1 \otimes x)(f_2 \otimes y) = f_1(.) f_2(\xinv .\,x) \otimes xy:
    z \mapsto f_1(z) f_2(\xinv zx)\: xy. \label{eq:muldub}
\ee
For the antipode we find
\be
S(f\otimes x) = f(x (.)^{-1} \xinv)\otimes \xinv : z \mapsto f(x \zinv \xinv)
\:\xinv.
\ee
The unit of ${\cal D}(G)$ is given by
\be
1\otimes e : z \mapsto e. \label{eq:undub}
\ee
The counit $\ep$ of ${\cal D}(G)$ is
\be
\ep(f\otimes x) = f(e).
\ee
For the R-matrix, which is an element of ${\cal D}(G)\otimes {\cal D}(G)$, we
have
\be
R = \sum_{x \in G} (\dl_x \otimes e) \otimes (1 \otimes x) 
\label{eq:rmatrix}
\ee
where $\dl_x$ is the Kronecker delta on $x \in G$, and thus a (basis) element
of $C(G)$. Hence $R(y,z) = e\otimes y$, for $y,z \in G$.

$C(G)$ becomes a Hopf $*$-algebra with
\be
f^*(x) := \overline{f(x)}. 
\label{eq:star1}
\ee
The corresponding Hopf $*$-algebra structure on $\CC[G]$ is given by
$x^* := \xinv$. Now ${\cal D}(G)$ has a $*$-algebra structure such that
$C(G)$ and $\CC[G]$ are $*$-subalgebras:
\be
(f\otimes x)^* = ((f\otimes e)(1 \otimes x))^* = (1\otimes x^*)(f^* \otimes
e)= \overline{f(x.\,\xinv)} \otimes \xinv.
\ee
We verify that
\be
((f_1\otimes x)(f_2 \otimes y))^* = (f_2\otimes y)^*(f_1\otimes x)^*,
\label{eq:starprod}
\ee
so we get a $*$-algebra structure on ${\cal D}(G)$. 
\begin{rem}
{\rm It has been shown by Majid, \cite{Majbook}, Proposition 7.1.4, and 
\cite{Majid94}, that the quantum double of any Hopf $*$-algebra naturally 
becomes a Hopf $*$-algebra.}
\end{rem}
We now give a slight reformulation of our model (9), (13), (19) for
${\cal D}(G)$ as a $*$-algebra. The new model will suggest
how to generalize the definition of this $*$-algebra to the case where 
$G$ is a locally compact group.
Observe that there is a linear bijection 
\be
\DG = C(G) \otimes \CC[G] \Longleftrightarrow C(G\times G).
\label{eq:bijecti}
\ee
For this bijection we can write:
\bea
\fox &\mapsto& \left((y,z) \mapsto f(y) \dl_x(z)\right) \nn\\
\sum_{z\in G} F(.\, ,z)\otimes z & \leftarrow& F
\eea
Then $C(G\times G)$ is a $*$-algebra with multiplication
\be
(F_1 \bullet F_2)(x,y) = \sum_{z\in G} F_1(x,z)
F_2(z^{-1}xz, z^{-1}y),
\label{eq:prod}
\ee
and $*$-structure
\be
 F^*(x,y) = \overline{F(y^{-1}xy,\yinv)}.
\label{eq:star}
\ee
Essentially the same algebra would be obtained with some constant 
nonzero factor added on the right hand side of Eq.(\ref{eq:prod}).
\section{The case of a locally compact group G}
In the following, compact or locally compact spaces are always supposed to
be Hausdorff.
If $X$ is a compact space then $C(X)$ will denote the space of continuous
functions on $X$. If $X$ is locally compact space then $C_c(X)$ will denote
the space of continuous functions with compact support on $X$ and $C_0(X)$
the space of continuous functions f on $X$ forwhich vanish $\lim_{\xi\to
\infty} f(\xi)=0$.
In all cases, $\|.\|_\infty$ will denote the sup-norm.

Let $G$ be a compact group with Haar measure $dg$.
Then the most straightforward generalisation for $\DG$ is
${\cal D}(G):=C(G\times G)$ with $*$-algebra structure given by
\be
(F_1\bullet F_2)(x,y):=\int_G F_1(x,z)\,F_2(z^{-1}xz,z^{-1}y)\,dz.
\label{eq:prodcomp}
\ee
and (\ref{eq:star}).
One can indeed verify that the axioms of
a $*$-algebra are satisfied. Note that this algebra in general has no
unit element.
\begin{rem}
{\rm It would be attractive to construct the quantum double associated with
a non-finite compact group $G$ in the spirit of Majid (\cite{Majbook}, p.296).
So work with two Hopf algebras in duality
(for finite $G$ these were $C(G)$ and $\CC[G]$), and make their tensor 
product into a Hopf algebra in the style of the quantum double construction.

However,
it is a problem which algebras we should choose, and whether these should
be Hopf algebras in the algebraic sense or only in the topological sense.
A next problem, for representation theory, is the type of irreducible
$*$-representations to be considered: all algebraic representations or
only those which are moreover continuous in some sense?
Also, will different choices of algebras in the quantum double construction
give rise to the same class of representations?
We might start with $C(G)$, which is an algebra and which is a coalgebra only
in the topological sense, but the Hopf algebra of
trigonometric polynomials on $G$ or (if $G$ is a compact Lie group)
the algebra $C^{\infty}(G)$ will also be candidates.
For a (topological) Hopf algebra in duality with the first chosen algebra
there are also many choices. Whatever we may choose, it has some 
arbitrariness. 

Our decision to follow another approach, namely to generalize the algebra 
structure of $C(G\times G)$ (obtained from the quantum double construction 
for finite $G$)
to the case of non--finite $G$, is motivated because we can then make contact
with an existing representation theory of transformation group algebras.

There are several other approaches in literature to the definition of
quantum double which look more conceptual, but which will not be followed
in the present paper. (We thank dr. Klaas Landsman for bringing some of these 
references to our attention.) The first approach of Podl\'es and Woronowicz 
\cite{PodWon} defines the double group of a compact matrix quantum group.
The second approach of Baaj and Skandalis \cite{BaSka}, see also the earlier
paper by Skandalis \cite{Ska}, defines the quantum double of a so--called
Kac system. They remark that their construction is compatible with
both the Drinfel'd \cite{Drin} and the Podl\'es--Woronowicz \cite{PodWon} 
construction of quantum double. A third approach, of Bonneau \cite{Bon}
defines a topological quantum double. Finally, in \cite{Mueg} M\"uger gives
the Von Neumann double of a locally compact group. However, these papers do 
not give a classification of irreducible $*$-representations of the quantum 
double.}
\end{rem}

If $G$ is a unimodular locally compact group with Haar measure $dg$
and if we put ${\cal D}(G):=C_c(G\times G)$, then formulas
(\ref{eq:prodcomp}) and
(\ref{eq:star}) still define a $*$-algebra structure on ${\cal D}(G)$.

If $G$ is a locally compact group with left Haar measure $dg$ and
Haar modulus $\Delta$ defined by $d(gh)=\Delta(h)\,dg$, then
${\cal D}(G):=C_c(G\times G)$ becomes a $*$-algebra with multiplication
(\ref{eq:prodcomp}) and $*$-structure given by
\be
F^*(x,y)=\Delta(y^{-1})\,\overline{F(y^{-1}xy,y^{-1})}.
\label{eq:locstar}
\ee
In fact, this last case is still a special case of a more general
$*$-algebra considered in literature: a {\sl transformation group algebra}.
Such an algebra is defined as follows:
\begin{defn}
Let $G$ be a locally compact group acting continuously on a locally
compact space $X$. Denote the action by
\be
(g, \xi) \mapsto g\xi : G\times X \to X, \qquad (g \in G,\, \xi \in X)
\ee
Then  $C_c(X\times G)$ is called a {\sl transformation group algebra} if
it is equipped with a multiplication and $*$-operation given by
\bea
(F_1 \bullet F_2)(\xi,y) = \int_{G} F_1(\xi,z) F_2(z^{-1}\xi ,z^{-1}y) dz \nn\\
F^*(\xi,y) = \overline{F(\yinv \xi,\yinv)} \Dl (\yinv)
\label{eq:operations}
\eea
\end{defn}
Straightfoward computations show that $C_c(X\times G)$ becomes a $*$-algebra.
We recover our earlier case ${\cal D}(G)$ when we take $X:=G$ and
$G$ acting on itself by conjugation.

These algebras first occur in literature in Dixmier \cite{Dix}. According
to Dixmier, the unimodular case was considered earlier in unpublished
work by Godement. The study of these algebras was continued by
Glimm \cite{gli62}. Afterwards, they were considered by many
authors, often in a more generalized form. See the survey paper
by Packer \cite{Pack} for references.
\begin{rem}
{\rm Suppose that $G$ is a finite group acting on a 
finite set $X$. Then we can take the set
\be
B := \{ \dl_{\eta} \otimes \dl_g \: |\: \eta \in X, g \in G \}
\ee
to be a basis of $C_c(X\times G)$ (which is here the space of all 
complex--valued functions on $X\times G$). In terms of these basis
elements, the operations in Eq.(\ref{eq:operations}) become
\bea
(\dl_{\eta}\otimes \dl_g)\bullet (\dl_{\eta'}\otimes \dl_{g'}) =
\dl_{\eta, g\eta'}\, (\dl_{\eta}\otimes \dl_{gg'}), \nn \\
(\dl_{\eta}\otimes \dl_g)^* = \dl_{g^{-1}\eta} \otimes \dl_{g^{-1}}.
\eea
There is also a unit element: $1 = \sum_{\xi \in X} \dl_{\xi}\otimes \dl_e$. 
Now it follows easily that the $\dubbelZ$--span of $B$, considered as an 
associative ring with unit and with involution, is a {\sl based ring}, as
defined by Lusztig \cite{Lus87}. For $\DG$ ($G$ finite group) this was
already observed in \cite{DPR}. }
\end{rem}
\begin{rem}
{\rm Suppose that a unimodular locally compact group 
$G$ continuously 
acts on a locally compact space $X$ and that $X$ is equipped with a 
$G$--invariant positive Borel measure. Then the algebra $C_c(X\times G)$ 
can be made into a pre--Hilbert space with inner product
\be
\langle F_1, F_2 \rangle := \int_X (F_1\bullet F_2^*)(\xi, e) d\xi = 
\int_G \int_X F_1(\xi,z) \overline{F_2(\xi,z)} d\xi dz,
\ee
and we can easily derive that
\be
\langle F F_1, F_2\rangle = \langle F_1, F^* F_2\rangle, \qquad (F, F_1,
F_2 \in C_c(X\times G)). 
\ee
In particular, if $G$ and $X$ are finite, then $C_c(X\times G)$ becomes a
semisimple algebra. (This is true for any based ring with finite
$\dubbelZ$--basis, see Lusztig \cite{Lus87}).}
\end{rem}
\section{Representation theory of transformation group algebras}
In the following a lcsc group (resp.\ space) will mean a topological
group (resp.\ space) which is
locally compact, Hausdorff and second countable. Also, Hilbert spaces
will be assumed to be separable.
These assumptions are for convenience. We have made no effort to
check to which extent the results below remain true without these assumptions.

Let a lcsc group $G$ act continuously on a lcsc space $X$ and consider the
transformation group algebra $C_c(X\times G)$ as above.
Glimm \cite{gli62} defines a norm on $C_c(X\times G)$ by
\be
\| F\|_1 := \int_G \| F(.\, ,z)\|_{\infty} \; dz
\label{eq:norm}
\ee
with $\|.\,\|_{\infty}$ the sup-norm on $X$. Then
\be
\| F_1 \bullet F_2 \|_1 \leq \| F_1\|_1 \| F_2\|_1 
\ee

We will now classify the irreducible $\|.\|_1$-bounded $*$-representations
of $C_c(X\times G)$ (up to equivalence)
under a certain assumption about the $G$-space $X$, namely that $X$ is 
{\sl countably separated} (see Definition \ref{thm:defcounsep}). 
This classification is due to Glimm \cite{gli62} Theorem 2.2,
(see the proof of his implication (2)$\Rightarrow$(3)
and take $K$ equal to ${\cal K}$).
As we will sketch below, the classification will follow from
Mackey's \cite{Mack49}, \cite{Mack58} imprimitivity theorem, 
see also {\O}rsted's
\cite{Oer}  short proof of Mackey's theorem. In fact, if $X$ is a transitive
$G$-space, then the classification result is equivalent to the imprimitivity
theorem.

In the following we will work with representations $\pi$
of various kinds of structured sets (a group, a $*$-algebra,
Borel algebra, etc.) on a Hilbert space $\cal H$.
In all cases under consideration
this will give rise to a $*$-closed subset $\{\pi(A)\}$ of the space
${\cal L}({\cal H})$ of all bounded linear operators on $\cal H$.
Then the commutant $R(\pi)$ and bicommutant $R(\pi)'$
of this set of operators are both a
{\sl Von Neumann algebra} (a weakly closed $*$-subalgebra of the
algebra ${\cal L}({\cal H})$). Note that $R(\pi)''=R(\pi)$, so
$R(\pi)$ and $R(\pi)'$ are the commutant of each other.
$R(\pi)$ is called the
{\sl Von Neumann algebra associated with the representation $\pi$}.
Representations can be classified according to their corresponding
Von Neumann algebras.

A representation $\pi$ is called {\sl irreducible} if $\{0\}$ and $\cal H$
are the only closed subspaces of $\cal H$ which are invariant under all
operators $\pi(A)$. A well known theorem says that $\pi$ is irreducible iff
$R(\pi)=\CC\,I$, iff $R(\pi)'={\cal L}({\cal H})$, see for instance
Theorem 4.7, Chap.VII in \cite{Mee79}.
\begin{defn}
A mapping 
\be
P:E \mapsto P_E : \{\mbox{Borel subsets of $X$}\}
\to \{\mbox{projection operators on Hilbert space ${\cal H}$}\}
\ee
is called a {\sl projection valued measure} if
\begin{itemize}
\item[i.]  $P_{\emptyset}=0, \; P_X = I$
\item[ii.] $ P_{E\cap F} = P_E P_F$
\item[iii.] $ P_E = \sum_{i=1}^{\infty} P_{E_i} \;\mbox{(strong 
convergence)} \;\mbox{if}\; E= \cup_{i=1}^{\infty}E_i \;
\mbox{(disjoint union)}$
\end{itemize}
If $v,w \in {\cal H}$, then $E\mapsto \langle P_E v,w\rangle$ is a complex 
measure on $X$, which we write as $dP_{v,w}(\xi)$. 
\end{defn}
The projection valued measures $P$ on $X$ are in one--to--one correspondence
with the non--degenerate $*$-representations $\pi$ of $C_0(X)$:
\be
\pi(f)=\int_X f(\xi)\,dP(\xi), \qquad (f\in C_0(X)),
\label{eq:piandP}
\ee
with the following interpretation:
\be
\langle \pi (f) v,w\rangle = \int_X f(\xi) dP_{v,w}(\xi),\qquad (v,w \in 
{\cal H}).
\label{eq:pivw}
\ee
The operators $P_E$ and the operators $\pi(f)$ generate the same 
Von Neumann algebra. 
\begin{defn}
Let $(G,X)$ be as before, $\tau$ be a unitary representation
of $G$ on ${\cal H}$, then a {\sl system of imprimitivity} (s.o.i.) for $\tau$
based on $X$ is given by $(G,X,P)$, where $P$ is a projection valued measure
on $X$ acting on ${\cal H}$ such that
\be
\tau(y) P_E \tau(y)^{-1} = P_{yE}, \qquad \forall y \in G,\; \forall\; \mbox{
Borel subsets}\; E \subset X.
\ee
We now say that $(\tau, P)$ furnishes a representation of $(G,X)$. 
\end{defn}
\begin{thm}[Glimm]\footnote{See \cite{gli62}, Theorem 1.51.} 
There is a one--to--one correspondence between
$*$-representations $\tau_0$ of $C_c(X\times G)$ which are bounded in the norm 
(\ref{eq:norm}) 
and representations $(\tau,P)$ of $(G,X)$ (both on the same
Hilbert space $\cal H$), in the following way
\be
\tau_0 (F) := \int_G \int_X  F(\xi,z)\;dP(\xi)\;\tau(z)\; dz
\label{eq:taunot}
\ee
which  has to be interpreted as
\be
\langle \tau_0 (F) v, w \rangle = \int_G \int_X \; F(\xi,z)\;dP_{\tau(z)v,w}
(\xi) \; dz,\quad \forall v,w \in{\cal H}.
\label{eq:inprod}
\ee
All such representations $\tau_0$ are norm-decreasing, i.e.
$\|\tau_0(F)\| \le \|F\|_1$.
The operators $\tau_0(F)$ and the operators $\tau(z)$ and $P_E$ generate
the same Von Neumann algebra. In particular, $(\tau,P)$ is irreducible iff
$\tau_0$ is irreducible.
An isometry of Hilbert spaces implements an equivalence between representations
$(\tau,P)$ and $(\tau',P')$ of $(X,G)$ iff it implements an equivalence of
the corresponding representations $\tau_0$ and $\tau_0'$ of 
$C_c(X\times G)$.
\end{thm}
\begin{rem}
{\rm In view of the one--to--one correspondence $\pi \leftrightarrow P$ given 
by Eq.(\ref{eq:piandP}), we can describe representations $(\tau,P)$ of
$(G,X)$ equivalently as so--called {\sl covariant representations} $(\tau,\pi)$
of $(G, C_0(X))$, where $\tau$ is a unitary representation of $G$ on
${\cal H}$ and $\pi$ is a nondegenerate $*$-representation of $C_0(X)$
on ${\cal H}$ such that
\be
\tau(y) \pi(f) \tau(y)^{-1}  = \pi(f(\yinv .)), \qquad \forall y \in G,
\forall f \in C_0(X).
\ee
More generally, Takesaki \cite{Tak67}, Def.3.1, defined covariant 
representations of $(G,A)$ where $A$ is a (possibly noncommutative)
$C^*$--algebra, and $G$ is a locally compact automorphism group of $A$.

Combining Eqs.(\ref{eq:pivw}) and (\ref{eq:inprod}) shows how a 
representation $\tau_0$ of $C_c(X\times G)$ can be obtained from
the corresponding covariant representation $(\tau, \pi)$ of $(G,C_0(X))$:
\be
\tau_0 (F) = \int_G \pi(F(.,z))\, \tau(z)\, dz
\label{eq:tauandpi}
\ee
which can also be interpreted weakly, like in Eq.(\ref{eq:inprod}).}
\end{rem}
Note that if $(\tau, P)$ is an irreducible representation of $(G,X)$, and
if $E, E'$ are $G$-invariant Borel sets, then $P_E$ commutes with $P_{E'}$
and with all $\tau(z)$, and we conclude that $P_E = \mbox{const.}I$,
{\em i.e.} $P_E = I$ or $P_E = 0$ ($P_E$ is a projection operator).
A projection valued
measure $P$ such that, for all $G$-invariant Borel sets $E$, $P_E=0$ or 
$I$, is called
{\sl ergodic}.
\begin{defn} 
$(G,X)$ as above, is called {\sl countably separated}, if
there are countably many Borel sets $B_1, B_2,\ldots\;$ in $X$ which are 
$G$-invariant, such that for every $G$-orbit ${\cal O}$ in $X$ we have that 
\be
{\cal O} = \cap_{{\cal O} \subset B_i} B_i.
\ee
\label{thm:defcounsep}
\end{defn}
This holds for instance, if $G$ and $X$ are compact and second countable. 
\vspace*{.2cm}\\
With $(G,X)$ as above, the following conditions are equivalent (see
Glimm \cite{gli61}, Theorem 1):
\begin{itemize}
\item[i.]
$(G,X)$ is countably separated;
\item[ii.]
The orbit space $G\backslash X$ is $T_0$ in the quotient topology. (A
topological space is $T_0$ if for any two distinct points at least one of
the points has a neighbourhood to which the other point does not belong.)
\item[iii.]
Each orbit in $X$ is relatively open in its closure.
\item[iv.]
For each $\xi\in X$ the map $zG_\xi\mapsto z\xi\:\colon\: G/G_\xi\to G\xi$
is a homeomorphism, where $G\xi$ has the relative topology of a subspace of
$X$. (Here $G_\xi$ denotes the stabilizer of $\xi$ in $G$.)
\end{itemize}
Hence, if $G$ is compact then $(G,X)$ is countably separated.
Also note that each $G$-orbit in $X$ is necessarily a Borel set and that
property (iv) above implies that the mapping $zG_\xi\mapsto z\xi$ is
a Borel isomorphism.
\begin{lem}\footnote{See for instance Van der Meer in \cite{Mee79}, Ch. XI, 
Lemma 3.3} 
If $(G,X)$ is countably separated and $P$ is an ergodic system of imprimitivity
on $\cal H$, then there is a unique $G$-orbit $\cal O$ in $X$ such
that $P_{\cal O}=I$ and $P_E=0$ if $E$ is the complement of $\cal O$.
(Then we say that $P$ is concentrated on $\cal O$.)
In particular, the conclusion holds if $(\tau,P)$ is an irreducible
representation of $(G,X)$.
\end{lem}
{}From now on assume that $(G,X)$ is countably separated.
Let $\cal O$ be a $G$-orbit in $X$. Take some $\xi\in {\cal O}$ and let
$G_\xi$ be the stabilizer of $\xi$ in $G$.
There is a one--to--one correspondence between representations
$(\tau,P)$ of $(G,X)$ concentrated on $\cal O$,
representations $(\tau,P')$ of $(G,{\cal O})$, and
representations $(\tau,P'')$ of $(G,G/G_\xi)$.
Here $P$ is related to $P'$ by
$P_E=P'_{{\cal O}\cap E}$ ($E$ Borel set of $X$),
and $P'$ is related to $P''$ via the homeomorphism, hence Borel isomorphism
$zG_\xi\mapsto z\xi$.
The three representations $(\tau,P)$, $(\tau,P')$ and $(\tau,P'')$
are associated with the same Von Neumann algebra.

Under this correspondence, equivalent representations of $(G,G/G_\xi)$
will give rise to equivalent representations of $(G,X)$ which are
concentrated on $\cal O$.
Conversely, if $(\tau,P)$ and $(\sigma,Q)$ are equivalent representation of
$(G,X)$ and if $P$ is concentrated on $\cal O$, then $Q$ will also be
concentrated on $\cal O$ and, under the above correspondence, the two
equivalent representations of $(G,X)$ will correspond to two
equivalent representations of $(G,G/G_\xi)$. 

Thus, the classification of irreducible $\|.\|_1$-bounded $*$-reps of
$C_c(X\times G)$ (up to equivalence) is reduced in several steps to the
problem of classifying the irreducible representations of $(G,G/H)$, where $H$
is a closed subgroup of $G$. Now we make contact with the notion of induced
representation of a locally compact group and with the imprimitivity
theorem.
\subsection{Connection with induced representations}
Let $\alpha$ be a unitary representation of $H$ on the Hilbert space 
$V_\alpha$.
Choose a  non-zero quasi-invariant measure $d\mu$ on $G/H$.
Then $d\mu(z\xi)=R(\xi,z)\,d\mu(\xi)$ ($z\in G)$, with $R$ a strictly
positive continuous function on $G/H\times G$ satisfying
\be
R(\xi,yz)=R(z\xi,y)\,R(\xi,z).
\label{eq:funcR}
\ee
If there is an invariant measure $\mu$ on $G/H$, which is certainly the
case if $H$ is compact, we can take $R = 1$. 
Introduce the following space of functions:
\bea
{\cal L}^2_{\al} (G, V_{\al}) := \{ f:G\to V_{\al} \;|\; f(gh) = \al (h^{-1})
f(g),\;\;\forall h\in H,\; \mbox{for almost all}\; g \in G,\nn\\ 
\mbox{and} \quad
\|f\|^2:= \int_{G/H} \|f(z)\|_{V_\alpha}^2\,d\mu(zH)<\infty\}
\eea
It has a positive semi-definite inner product
\be
\langle f_1,f_2\rangle:=
\int_{G/H}\langle f_1(z),f_2(z)\rangle_{V_\alpha}\,d\mu(zH).
\ee
We obtain a Hilbert space by taking the quotient space with respect
to the subspace of functions with norm zero:
\be
L^2_{\al} (G,V_{\al}) = {\cal L}^2_{\al} (G, V_{\al})\: /\: \{f \in 
{\cal L}^2_{\al} (G, V_{\al})\:|\: \|f\| = 0\}
\ee
Now we can write down the following unitary representation of $G$ on 
$L^2_{\al}(G,V_{\al})$:
\be
(\tau_\alpha(y)f)(x):=(R(xH,y^{-1}))^{1/2}\,f(y^{-1}x),
\ee
and the following projection valued measure:
\be
(P^{\al}_E f)(x) = \chi_E(xH) f(x), \qquad{E\;\mbox{Borel subset in}\; G/H}
\label{eq:projmeas}
\ee
where $\chi_E$ is the characteristic function of $E$ (namely, $\chi_E(x)=1$
if $x \in E$, and zero elsewhere). 
\begin{prop}
$(\tau_{\al},P^{\al})$ is a unitary representation of $(G,G/H)$ on 
$L^2_{\al}(G,V_{\al})$. The Von Neumann algebras $R(\al)$ and 
$R(\tau_\alpha,P^\alpha)$ are isomorphic. In particular, $\al$ is 
irreducible iff $(\tau_\alpha,P^\alpha)$ is irreducible.
Also, the equivalence class of $\al$ corresponds to the equivalence class
of $(\tau_\alpha,P^\alpha)$.
\end{prop}
Now we can use Mackey's imprimitivity theorem, adapted to our specific
situation:
\begin{thm}[Mackey]  
If $(\tau, P)$ is a representation of $(G, G/H)$, then $(\tau, P)$ is 
equivalent to an (induced) representation
$(\tau_{\al}, P^{\al})$, with $\al$ a unitary representation of $H$. 
$(\tau, P)$ is irreducible iff $\al$ is irreducible.
\end{thm}
Summarizing we have the following equivalences:\\
Irreducible representation (irrep) $\tau_0$ of $C_c(X\times G)$ 
$\Longleftrightarrow$ \\
irrep $(\tau, P)$ of $(G,X)$ $\Longleftrightarrow$ (if $(G,X)$ is 
countably separated) \\
irrep $(\tau, P)$ of $(G,{\cal O})$, with ${\cal O} \simeq G/H$, and
$H$ the stabilizer of point $\xi_0 \in X  \Longleftrightarrow $ 
(imprimitivity)\\
irrep $(\tau_{\al}, P^{\al})$ of $(G,{\cal O})$, with $\al$ unitary
irrep of $H$. 
\subsection{Induced representations of transformation group algebras}
Let $\tau_{\al,0}$ be the representation of $C_c(X\times G)$ obtained by
extending the representation $(\tau_{\al},P^{\al})$ of $(G,{\cal O})$ to
$(G,X)$ (putting $P^{\al}=0$ on the complement of ${\cal O}$) and next
lifting it to a representation of $C_c(X\times G)$. Then it follows from 
(\ref{eq:inprod}) that
\be
\langle \tau_{\alpha,0}(F) \phi,\psi\rangle=
\int_G\int_X F(\xi,z)\,dP_{\tau_\alpha(z)\phi,\psi}^\alpha(\xi)\,dz,
\label{eq:inpphipsi}
\ee
where $F\in C_c(X\times G)$, $\phi,\psi\in L_\alpha^2(G,V_\alpha)$.
Now use that
\be
\int_X f(\xi)\,dP_{\phi,\psi}^\alpha(\xi)=
\int_{G/H} f(x\xi_0)\,\langle \phi(x),\psi(x)\rangle\,d\mu(xH),
\label{eq:intf}
\ee
where $f\in C_c(X)$, $\phi,\psi\in L_\alpha^2(G,V_\alpha)$, $\xi_0 \in X$.
This follows since, for a Borel set $E$ in $X$,
\bea
\begin{array}{lll}
P_{\phi,\psi}^\al(E)=&\langle P_E^\alpha\phi,\psi\rangle= \nn \\
& \left\{\begin{array}{l}
0 \qquad \qquad \mbox{if $E$ is in the complement of $\cal O$},\nn\\
\int_{G/H}\chi_E(xH)\,\langle\phi(x),\psi(x)\rangle\,d\mu(xH)=
\int_E\langle\phi(x),\psi(x)\rangle\,d\mu(xH)\nn\\
\qquad   \mbox{if $E$ is a Borel set of $\cal O$ and is transfered 
to a Borel set of $G/H$.}
\end{array} \right.\end{array}
\eea
Here we used formula (\ref{eq:projmeas}). From (\ref{eq:inpphipsi}) and
(\ref{eq:intf}), we obtain:
\be
\langle \tau_{\alpha,0}(F) \phi,\psi\rangle=
\int_G\int_{G/H} F(x\xi_0,z)\,\langle(\tau_\alpha(z)\phi)(x),\psi(x)\rangle\,
d\mu(xH)\,dz.
\ee
So finally we arrive at
\be
(\tau_{\al,0}(F) \phi)(x) :=\int_G F(x\xi_0,z) (\tau_\alpha(z)\phi)(x)\:dz.
\ee
\begin{thm}  
Let $(G,X)$ be countably separated. Let $\{{\cal O}_A\}_{A\in{\cal A}}$
be the collection of $G$-orbits in $X$ ($\cal A$ an index set).
For each $A\in{\cal A}$ choose some $\xi_A\in{\cal O}_A$,
let $N_A$ be the stablizer of $\xi_A$ in $G$, choose some quasi-invariant
measure $\mu_A$ on $G/N_A \simeq {\cal O}_A$ and let $R_A$ be the 
corresponding $R$-function given by (\ref{eq:funcR}).
For each $\alpha\in\widehat{N_A}$ choose a representative, also
denoted by $\alpha$, which is an irreducible unitary representation of 
$N_A$ on some Hilbert space $V_\alpha$.
Then, for $A\in{\cal A}$ and $\alpha\in\widehat{N_A}$ we have
mutually inequivalent irreducible $\|.\|_1$-bounded $*$-representations
$\tau_\alpha^A$ of $C_c(X\times G)$ on $L_\alpha^2(G,V_\alpha)$ given by
\be
(\tau_\alpha^A(F)\phi)(x):=
\int_G F(x\xi_A,z)\,(R_A(x \xi_A,z^{-1}))^{1/2}\,\phi(z^{-1}x)\,dz,
\label{eq:ourrep}
\ee
and all irreducible $\|.\|_1$-bounded $*$-representations of $C_c(X\times G)$
are equivalent to some $\tau_\alpha^A$.
\label{thm:rep}
\end{thm}
Proof: follows from the statements before.
\vspace*{.1cm}\\
For purposes of later reference we formulate the specialization of
Theorem \ref{thm:rep} to the case
${\cal D}(G):=C(G\times G)$ ($G$ compact group)
with $*$-algebra operations given by (\ref{eq:prodcomp}) and (\ref{eq:star}), 
and with norm $\|.\|_1$ defined by (\ref{eq:norm}). Since $G$ is compact, the
requirement of countable separability is certainly fulfilled. The orbits
${\cal O}_A$ are identified as the conjugacy classes $C_A$. 
\begin{cor}
Let $\{C_A\}_{A\in{\cal A}}$
be the collection of conjugacy classes of the compact group
$G$ ($\cal A$ an index set). For each $A\in{\cal A}$ choose some 
$g_A\in C_A$ and let $N_A$ be the centralizer of $g_A$ in $G$.
For each $\alpha\in\widehat{N_A}$ choose a representative, also
denoted by $\alpha$, which is an irreducible unitary representation of
$N_A$ on some Hilbert space $V_\alpha$.
Then, for $A\in{\cal A}$ and $\alpha\in\widehat{N_A}$ we have
mutually inequivalent irreducible $\|.\|_1$-bounded $*$-representations
$\tau_\alpha^A$ of $C(G\times G)$ on $L_\alpha^2(G,V_\alpha)$ given by
\be
(\tau_\alpha^A(F)\phi)(x):=
\int_G F(xg_A x^{-1},z)\,\phi(z^{-1}x)\,dz,
\label{eq:repcomp}
\ee
and all irreducible $\|.\|_1$-bounded $*$-representations of $C(G\times G)$
are equivalent to some $\tau_\alpha^A$.
\label{thm:cor}
\end{cor}
\begin{rem}
{\rm Because of Eq.(\ref{eq:tauandpi}), we can describe the representation 
$\tau^A_{\al}$ in Eq.(\ref{eq:ourrep}) equivalently as the covariant
representation $(\tau_{\al}, \pi^A_{\al})$ of $(G, C_0(X))$, where
$\tau_{\al}$ is the unitary representation of $G$ which is induced by
the representation $\al$ of $N_A$, and where
\be
(\pi^A_{\al}(f)\phi)(x) := f(x \xi_A) \phi(x), \qquad (\phi \in 
L^2_{\al}(G, V_{\al}), \, f\in C_0(X)). 
\ee
By the definition of induced covariant representations in \cite{Tak67}, the
covariant representation $(\tau_{\al}, \pi^A_{\al})$ of $(G, C_0(X))$ is
induced by the covariant representation $(\al, \xi_A)$ of $(N_A, C_0(X))$
on $V_{\al}$, where} 
\be
\xi_A (f) v := f(\xi_A) v, \qquad (v\in V_{\al},\, f\in C_0(X)).
\ee
\end{rem}
\begin{rem} 
{\rm In the situation of Theorem \ref{thm:rep}, a {\sl Borel cross--section} 
for $G/N_A \simeq {\cal O}_A$
means a Borel mapping $s_A:{\cal O}_A \rightarrow G$ such that $s_A(\xi) 
\xi_A = \xi$, for all $\xi \in {\cal O}_A$. By a theorem of Mackey (see
\cite{Mack76}), such a Borel cross section always exists. In terms of
$s_A$ the representation $\tau^A_{\al}$ in (\ref{eq:ourrep}) can be
equivalently described as a $*$-representation acting on $L^2({\cal O}_A,
V_{\al};\mu_A)$ (where $\mu_A$ is an invariant measure on the orbit 
(or conjugacy class) ${\cal O}_A$) by
\bea
(\tau^A_{\al} (F) \phi)(\xi) = \int_G F(\xi,z) R_A(\xi,z^{-1}) 
\al(s_A(\xi)^{-1}z\,s_A(z^{-1}\xi))\, \phi(z^{-1}\xi)\, dz \nn\\
 \xi\in  {\cal O}_A, \phi\in L^2({\cal O}_A, V_{\al};\mu_A).
\label{eq:sectrep}
\eea
Consider in particular the case of $A \in {\cal A}$ such that ${\cal O}_A$ 
has only one element $\xi_A$ (equivalently, $N_A = G$). Then $R_A$ is 
identically 1, we can take $s_A(\xi_A) = e$, and we see from 
Eq.(\ref{eq:sectrep}) that $\tau^A_{\al}$ is a $*$-representation of
$C_c(G\times G)$ on $V_{\al}$ given by}
\be
\tau^A_{\al} (F) = \int_G F(\xi_A,z) \al(z)\,dz,\quad N_A=G,\;\al\in\widehat{G}.
\ee
\label{thm:specrep}
\end{rem}
\section{The case of finite $G$}
As an example we specialise to the case of finite group $G$, where $G$ acts
on itself via conjugation. We derive the 
unitary irreducible representations in the way outlined in the last section,
and will show that our result is isomorphic to the representations derived 
in \cite{DPR}.

For a finite group the space $X = G$ is countably separated. 
The Hilbert space is
\be
\HAal := \{ v \in L^2(G;V_{\al})\,|\,\forall n \in N_A: v(xn) = \al(\ninv)
v(x)\,\, \mbox{for almost all $x$} \}. 
\label{eq:hilbwe}
\ee
Corollary \ref{thm:cor} holds, and Eq.(\ref{eq:repcomp}) can now be
rewritten as
\be
(\pi^A_{\al}(F)v)(x)=\sum_{y\in G} F(x\age \xinv,y) v(\yinv x), 
\qquad F\in C(G\times G),
\label{eq:finrepwe}
\ee
where $v$ is in the Hilbert space (\ref{eq:hilbwe}).
In \cite{DPR} these representations were derived by means of induction of
algebra representations. For completeness we will show here how that works
in this particular example.

Let
\be
{\cal B}_A := C(G) \otimes \CC[N_A],
\ee
considered as a subalgebra of $\DG$. Denote a general element of a spanning 
set for this subalgebra by $f\otimes n$, with $f \in C(G)$ and $n \in N_A$. 
Then define the representation of ${\cal B}_A$ on $V_{\al}$:
\be
\Pi_{\al} (f \otimes n) v := f(\age) \al(n) v,\;\;\; v\in V_{\al}.
\ee
This is indeed a representation, since
\bea
\Pi_{\al}(f_1 \otimes n_1) \Pi_{\al}(f_2 \otimes n_2) = f_1(\age) f_2(\age)
\al(n_1 n_2) &=& \Pi_{\al}(f_1(.) f_2(n_1^{-1}.\,n_1) \otimes n_1 n_2) \nn\\
&=& \Pi_{\al}((f_1 \otimes n_1)(f_2 \otimes n_2)).
\eea
Now we induce this representation of ${\cal B}_A$ on $V_{\al}$ to a 
representation $\PAal$ of $\DG$ on the representation space
\be
V^A_{\al} := \DG \otimes_{{\cal B}_A} V_{\al},
\label{eq:vectsp}
\ee
that is, $V^A_{\al}$ is a left module of $\DG$ 
\footnote{For the case of $G$ a finite group this tensor product is well 
defined. This procedure is also valid for
the case where $G$ is a locally compact group. Then by this representation 
space we mean a certain well--chosen completion of the
tensor product space.}. 
We note that a general element of (a spanning set of) this
representation space can be written in the following way:
\be
(\fox)\otimes_{{\cal B}_A} v = (1\otimes x)(f(x.\,x^{-1})\otimes e)
\oba v = f(x\age x^{-1}) (1 \otimes x) \oba v
\ee
This is effectively an element of $(1\otimes \CC[G])\oba V_{\al}$, which 
equals $\CC[G]\otimes_{\al} V_{\al}$, where $\otimes_{\al}$ denotes the
tensor product with equivalence relation
\be
xh\otimes v \equiv x\otimes \al(h) v, \quad h\in N_A.
\ee
Therefore, there is a bijection
\bea
\DG \oba V_{\al} &\Longleftrightarrow& \CC[G] \otimes_{\al} V_{\al} \\
(\fox) \oba v &\mapsto& f(x\age \xinv) x \otimes_{\al} v \nn \\
(1\otimes x) \oba v &\leftarrow& x \otimes_{\al} v \nn
\label{eq:biject}
\eea
Transfer the representation $\PAal$ from $\VAal$ to $\CC[G]\otimes_{\al}
V_{\al}$ under this bijection, and call the resulting representation
again $\PAal$. Because
\be
\PAal (1\otimes y)(x\otimes_{\al} v) = yx \otimes_{\al} v
\label{eq:grouprep}
\ee
$\PAal (1\otimes .)$ is the representation of $G$ induced by the 
representation $\al$ of $N_A$ on $V_{\al}$. For $C(G)$ we have 
\be
\PAal (f\otimes e) (x \otimes_{\al} v) = f(x\age \xinv) x\otimes_{\al} v.
\label{eq:finrepthey}
\ee
and so the representation $\PAal$ can be completely described by inducing
the representation $\al$ of $N_A$ to $G$ and by specifying the representation
$\PAal(.\otimes e)$ of $C(G)$ on the representation space of the 
representation of G obtained by inducing $\al$. 

In order to show that the representation in Eqs.(\ref{eq:grouprep}), 
(\ref{eq:finrepthey}) is equivalent to the representation in 
Eq.(\ref{eq:finrepwe}), we give the linear intertwining bijections between 
the Hilbert spaces from Eqs.(\ref{eq:hilbwe}) and (\ref{eq:vectsp}). For 
$x \otimes_{\al} v \in \CC[G]\otimes_{\al} V_{\al}$, and $w \in \HAal$ 
as in Eq.(\ref{eq:hilbwe}) we have:
\bea
x\otimes_{\al} v &\mapsto& \frac{1}{N_A} \sum_{n\in N_A} \dl_{x n^{-1}}(.)
\al(n) v \nn \\
\sum_{x \in G} x\otimes_{\al} w(x) &\leftarrow& \quad w(.)
\eea
The last mapping shows that if $\bar{x}$ is taken to be the representative 
of the right coset $xN_A$, and
$\{e_k\}_{k=1}^{dim V_{\al}}$ a basis of $V_{\al}$ with $v = v^k e_k$, then
$\bar{x} \otimes e_k$ is a basis for $V^A_{\al}$, and we can expand $v$ on
this basis:
\be
v = \sum_{\bar{x} \in G/N_A} \sum_k v^k(\bar{x})\; \bar{x} \otimes e_k
\ee
For completeness we mention that there are also two other vector spaces
which are isomorphic to $V^A_{\al}$. We define a mapping 
$s: G/N_A \to G$ such that $s(gN_A)N_A = gN_A$ for all $gN_A \in G/N_A$. Then
\begin{itemize}
\item[i.] $V^A_{\al,s} := \CC[G/N_A] \otimes V_{\al}$
\item[ii.] $\widetilde{V^A_{\al,s}} := \{\Phi: G/N_A \to V_{\al}\}$
\end{itemize}
and we can find compatible linear isomorphisms between them and the 
representation spaces $V^A_{\al}$ and $\HAal$. This concludes our discussion
of the connection with \cite{DPR}.
\section{Examples}
\subsection{$SU(2)$}
To illustrate the construction of the quantum double for a compact group $G$
and its irreducible representations, we now consider the case of $G = SU(2)$.

Let 
\bea
g_{\theta} = \left(\begin{array}{cc} e^{i\theta} & 0 \\ 0 & e^{-i\theta}
\end{array} \right) \in SU(2)
\eea
be the representative of the conjugacy class $C_{\theta} := \{ g g_{\theta}
g^{-1}\,|\, g\in SU(2)\}$, with $0\leq\theta\leq\pi.$  For $0<\theta<\pi$,
the centralizer $N_{\theta}$ equals $U(1)$, which is embedded in $SU(2)$ 
in the following way: $U(1) := \{ g_{\theta}\,|\, -\pi<\theta\leq \pi\}$.
For $\theta = 0, \pi$ the centralizers $N_0, N_{\pi}$ are equal to $SU(2)$. The
irreducible unitary representations of $U(1)$ are given by $\al_n(g_{\theta}) =
e^{in\theta}, (n\in\dubbelZ)$. Let
\be
L^2_n(SU(2)) := \{\phi\in L^2(SU(2))\;|\; \phi(g g_{\theta}) = 
\al_n(g_{\theta}^{-1}) \phi(g),\quad \forall \theta \in [0,2\pi]\}.
\ee
Then, by Corrolary \ref{thm:cor}, we have for $0<\theta<\pi$ and
$n\in\dubbelZ$ irreducible $*$-representations $\tau^{\theta}_n$ of
${\cal D}(SU(2))$ on $L^2_n(SU(2))$ given by
\be
(\tau^{\theta}_n(F)\phi)(x) = \int_{SU(2)} F(x g_{\theta}\xinv, z) 
\phi(z^{-1}x)\, dz, \qquad x\in SU(2), \phi\in L^2_n(SU(2)).
\ee
For $\theta = 0, \pi$ we find by Remark \ref{thm:specrep} the irreducible 
$*$-representations
\be
\tau^{\theta}_l(F) = \int_{SU(2)} F(g_{\theta},z) \pi_l(z)\,dz,
\qquad l=0,\half,1,\ldots 
\ee
on $\CC^{2l+1}$, where $\pi_l(z)$ denotes the $(2l+1)$-dimensional unitary 
irreducible representation of $SU(2)$.

Up to equivalence, the representations $\tau^{\theta}_n,$ with $(0<\theta<\pi,
\,n\in \dubbelZ)$ and $\tau^{\theta}_l, (\theta = 0, \pi,\, 
l=0,\half,1,...)$ give all irreducible $\|.\|_1$-bounded $*$-representations 
of ${\cal D}(SU(2))$.
\subsection{$SL(2,\dubbelR)$}
As an example of a non--compact group we take $SL(2,\dubbelR)$. The
(generalized) eigenvalues $\lm_1, \lm_2$ of elements of this group are either 
real, or complex conjugated, and $\lm_1 \lm_2 = 1$. We see that the space of
conjugacy classes can be split in the following way, according to the
eigenvalues of elements of the classes:
\begin{itemize}
\item[i.\& ii.]  For $\lm_1 = e^{i\theta},\; \lm_2 = e^{-i\theta} \;$
($0<\theta<\pi)$, we find two disjoint conjugacy classes:\\
$C_{\theta}:= \{ g\, u_{\theta}\,g^{-1}\;|\; g\in SL(2,\dubbelR)\}$, 
and $C_{-\theta}:=\{ g\, u_{-\theta}\,g^{-1}\;|\; g\in SL(2,\dubbelR)\}$, with
\bea
u_{\psi} = \left(\begin{array}{cc} \cos\psi & \sin\psi \\ 
-\sin\psi & \cos\psi
\end{array} \right) \in SL(2,\dubbelR)
\label{eq:utheta}
\eea
They have the same centralizer $N_{\theta} = N_{-\theta} = U(1) := 
\{ u_{\psi}\,|\, -\pi<\psi\leq \pi\}$, embedded in $SL(2,\dubbelR)$. 
This centralizer has irreducible unitary representations 
labeled by $n \in \dubbelZ$, like in the example of $SU(2)$. 
\item[iii.] For $\lm_1 = e^t, \lm_2 = e^{-t}, t>0$, there is one conjugacy
class:\\
$C_t:= \{ g\, a_t\, g^{-1}\;|\; g\in SL(2,\dubbelR)\}$, with
\bea
a_s = \left(\begin{array}{cc} e^s & 0 \\ 0 & e^{-s} \end{array} \right) 
\in SL(2,\dubbelR).
\label{eq:at}
\eea
This has centralizer $N_t := \{\pm a_s\;|\; s\in \dubbelR\}$, 
which means that $N_t\simeq \dubbelR\times\dubbelZ_2$. The irreducible
unitary representations of $N_t$ are labeled by the pairs $(b,\ep)$, with
$ b \in \dubbelR,\; \ep = \pm 1$.
\item[iv.] $\lm_1 = -e^t, \lm_2 = -e^{-t}$ is associated to the class \\
$\bar{C}_t := \{ g\,\bar{a}_t\, g^{-1}\;|\; g\in SL(2,\dubbelR)\}$ and
\bea
\bar{a}_s =\left(\begin{array}{cc} -e^s & 0 \\ 0 & -e^{-s} \end{array} 
\right) \in SL(2,\dubbelR)
\eea
It has the same centralizer as $C_t$.
\item[v. \& vi.]  In case the eigenvalues are $\lm_1 = \lm_2 = 1$ we 
distinguish two conjugacy classes $C_e$ and $C_1$:
$C_e := \{I\}$, which has centralizer $N_e := SL(2,\dubbelR)$. The unitary 
irreducible representations of $SL(2,\dubbelR)$ have been classified by
Bargmann \cite{Bar}, see for instance Van Dijk in \cite{Mee79}.\\
$C_1 := \{ g\, n_1\, g^{-1}\;|\; g\in SL(2,\dubbelR)\}$,
\bea
n_1 = \left(\begin{array}{cc} 1 & 1 \\ 0 & 1
\end{array} \right) \in SL(2,\dubbelR).
\label{eq:n1}
\eea
This class has centralizer $N_1$ consisting of matrices
\bea
\pm \left(\begin{array}{cc} 1 & z \\ 0 & 1 \end{array} \right),
\qquad z\in \dubbelR
\eea
and thus $N_1\simeq \dubbelR \times \dubbelZ_2$. The irreducible unitary 
representations of $N_1$ are labeled by $(d,\ep)$, with $d\in \dubbelR, 
\ep = \pm 1$.
\item[vii. \& viii.] Finally for $\lm_1 = \lm_2 = -1$ we find the classes
$C_{-e}:= \{-I\}$, again with centralizer $SL(2,\dubbelR)$, and 
$\bar{C}_1 :=\{ g\,\bar{n}_1\, g^{-1}\;|\; g\in SL(2,\dubbelR)\}$,
\bea
\bar{n}_1 = \left(\begin{array}{cc} -1 & 1 \\ 0 & -1
\end{array} \right) \in SL(2,\dubbelR).
\eea
with centralizer $\bar{N}_1 = N_1$.
\end{itemize}
A straightforward argument from linear algebra shows that the various classes
given above are indeed distinct conjugacy classes, and that they are all
conjugacy classes.

Next we prove that $SL(2,\dubbelR)$ is countably separated, by using
the second equivalence under Definition \ref{thm:defcounsep}.  To that aim
we note that taking the trace $tr$ is a continuous function on the space
of conjugacy classes, so conjugacy classes with different values of $tr$
are lying in disjoint open subsets. There are two possible types of 
obstructions:

The trace equals $(2\cos \theta)$ on $C_{\theta}$ and $C_{-\theta}$, 
$(0<\theta<\pi)$. However, then the preimage for the continuous function
\bea 
\left(\begin{array}{cc} a & b \\ c & d \end{array}\right) \mapsto b
\eea
on $SL(2,\dubbelR)$ takes $b>0$ to an open subset containing  
$C_{\theta}$, and $b<0$ to an open subset containing $C_{-\theta}$.
So by considering inverse images $C_{\theta}$ and $C_{-\theta}$ are in
disjoint open subsets.

The trace equals 2 on $C_e$ and $C_1$. However, $C_1$ is included
in the open subset $SL(2,\dubbelR)\backslash \{I\}$ of $SL(2,\dubbelR)$. 
This means that in the quotient topology it lies in an open subset which 
does not include $C_e=\{I\}$. Similarly if $tr = -2.$\\
This means that the space of conjugacy classes is $T_0$, and thus that
$SL(2,\dubbelR)$ is countably separated. $\qquad \qquad \qquad\Box$
\vspace{.2cm}\\

This classification of conjugacy classes and the representations of their 
centralizers enables us to classify the irreducible $*$-representations
of ${\cal D}(SL(2,\dubbelR))$. Using Corrolary \ref{thm:cor} for the
cases (i) and (ii) listed above, we find for example
\be
(\tau^{\theta}_n(F)\phi)(x) = \int_{SL(2,\dubbelR)} F(x u_{\theta}\xinv, z) 
\phi(z^{-1}x)\, dz, \qquad x\in SL(2,\dubbelR), \phi\in L^2_n(SL(2,\dubbelR)),
\ee 
for $-\pi<\theta<\pi, \theta\neq 0$. On the various (nontrivial) conjugacy 
classes there is an invariant measure, and therefore the $R$-function is 
equal to 1. 

For $C_e$ and $C_{-e}$ it follows from Remark \ref{thm:specrep} that
\be
\tau^{\pm e}_r(F) = \int_{SL(2,\dubbelR)} F(\pm I,z) r(z)\,dz
\ee
where $r(z)$ denotes a unitary irreducible representation of $SL(2,\dubbelR)$.
\subsection*{Acknowledgements}
The second author was supported by the Dutch Science Foundation FOM/NWO.
We would like to thank professor Sander Bais for stimulating discussions
and putting us on the track of the subject of this paper
\end{document}